\documentclass[11pt]{amsart}
\usepackage{amsfonts}
\usepackage{amsfonts}
\usepackage[arrow,matrix]{xy}
\usepackage{amsmath,amssymb,bbm, amscd, amsthm,mathrsfs,hyperref}
\usepackage{longtable}
\theoremstyle{plain}
\newtheorem{thm}{Theorem}[section]
\newtheorem{cor}[thm]{Corollary}
\newtheorem{lem}[thm]{Lemma}

\theoremstyle{definition}
\newtheorem{defn}[thm]{Definition}

\newtheorem{rem}[thm]{Remark}

\def\al{\alpha}
\def\bt{\beta}
\def\dt{\delta}

\def\sg{\sigma}

\def\vph{\varphi}
\def\lmd{\lambda}

\def\vps{\varepsilon}
\def\om{\omega}

\def\bdt{\Delta}

\def\mc{\mathcal}

\def\ot{\otimes}

\def\se{\leqslant}
\def\le{\geqslant}

\def\ra{\rightarrow}
\def\lra{\longrightarrow}
\def\la{\leftarrow}
\def\xra{\xrightarrow}
\def\ced{\centerdot}
\def\Hom{\operatorname {Hom}}
\def\Ext{\operatorname {Ext}}
\def\Tor{\operatorname {Tor}}

\def\id{\operatorname {id}}

\def\Mod{\operatorname {Mod}}
\def\H{\operatorname {H}}

\def\GL{\operatorname {GL}}
\def\tr{\operatorname {tr}}
\def\ob{\operatorname {ob}}
\def\t{\text}

\def\it{\textit}

\def\kk{\mathbbm{k}}

\begin{document}
\title[Hopf-Galois objects of Calabi-Yau Hopf algebras]{\bf Hopf-Galois objects of Calabi-Yau Hopf algebras}

\author{Xiaolan YU}
\address {Xiaolan YU\newline Department of Mathematics, Hangzhou Normal University, Hangzhou, Zhejiang 310036, China}

\email{xlyu@hznu.edu.cn}

\date{}

\begin{abstract}
By using the language  of cogroupoids, we show that Hopf-Galois objects of a twisted Calabi-Yau Hopf algebra with bijective antipode are still twisted Calabi-Yau, and give their Nakayama automorphism explicitly. As applications, cleft Galois objects of twisted  Calabi-Yau Hopf algebras and Hopf-Galois objects of the quantum automorphism groups of non-degenerate bilinear forms are proved to be twisted Calabi-Yau.

\end{abstract}

\keywords{Hopf-Galois object; Calabi-Yau algebra; Quantum group}
\subjclass[2000]{16E65, 16E40, 16W30, 16W35.}

\maketitle

\section*{Introduction}
We work over a fixed field $\kk$. Let $E\in \t{GL}_m(\kk)$ with $m\le 2$ and let $\mc{B}(E)$ be the algebra  presented by generators $(u_{ij})_{1\se i,j\se m}$ and   relations
$$E^{-1}u^tEu=I_m=uE^{-1}u^tE,$$
where $u$ is the matrix $(u_{ij})_{1\se i,j\se m}$, $u^t$ is the transpose of $u$ and $I_m$ is the identity matrix. The algebra $\mc{B}(E)$ is a Hopf algebra and was defined by Dubois-Violette and Launer \cite{dvl} as the quantum automorphism group of the non-degenerate bilinear form associated to $E$. It is not difficult to check that $\mc{B}(E_q)=\mc{O}_q(\t{SL}_2(\kk))$ (the quantised coordinate algebra of $\t{SL}_2(\kk)$), where $$E_q=\left(\begin{array}{cc}0&1\\-q^{-1}&0\end{array}\right).$$
So the Hopf algebras $\mc{B}(E)$ are generalizations  of $\mc{O}_q(\t{SL}_2(\kk))$.

For any $E\in \t{GL}_m(\kk)$, in  \cite{bi}, Bichon constructed a free Yetter-Drinfeld resolution of the trivial module over $\mc{B}(E)$.   Consequently, $\mc{B}(E)$ is smooth of dimension 3 and satisfies a Poincar\'{e} duality between Hochschild homology and cohomology. Namely, $\mc{B}(E)$ is a twisted Calabi-Yau (CY for short) algebra of dimension 3 (Definition \ref{tcy}). CY algebras were introduced by Ginzburg in 2006 \cite{g2}. They are studied in recent years because of their applications in algebraic geometry and mathematical physics. Twisted CY algebras process similar homological properties as the CY algebras and include CY algebras as a subclass. They are the natural algebra analogues of the Bieri-Eckmann duality groups \cite{be}. Associated to a twisted CY algebra, there exists a so-called Nakayama automorphism. The Nakayama automorphism is the tool that enables one to describe the duality between Hochschild cohomology and homology. This automorphism is unique up to an inner automorphism.
A twisted CY algebra is CY if and only if its Nakayama automorphism is an inner automorphism.

Aubriot gave the classification of Hopf-Galois objects (Definition \ref{hg}) of the algebras $\mc{B}(E)$ in \cite{au}. Let $E\in \t{GL}_m(\kk)$ and $F\in \t{GL}_n(\kk)$ with $m,n\le2$. Define $\mc{B}(E,F)$ to be the algebra with generators $u_{ij}$, $1\se i\se m$, $1\se j\se n$, subject to the relations:
$$F^{-1}u^t E u=I_n;\;\;uF^{-1}u^t E=I_m.$$
On one hand, if $\tr(E^{-1}E^t)=\tr(F^{-1}F^t)$, then $\mc{B}(E,F)$  is a left $\mc{B}(E)$-Galois object. On the other hand, if $A$ is a left $\mc{B}(E)$-Galois object, then there exists $F\in \t{GL}_n(\kk)$ satisfying $\tr(E^{-1}E^t)=\tr(F^{-1}F^t)$, such that $A\cong \mc{B}(E,F)$ as Galois objects. The algebras $\mc{B}(E)$ are twisted CY algebras. One naturally asks whether $\mc{B}(E,F)$ are twisted CY as well. More generally, the CY property of Hopf algebras was discussed in \cite{bz} and \cite{hvz},  we would like to know whether a Hopf-Galois object of a twisted CY Hopf algebra is still a twisted CY algebra.

Hopf-Galois objects can be described by the language of cogroupoids. A cogroupoid $\mc{C}$ consists of a set of objects $\ob(\mc{C})$ such that for any $X,Y\in\ob(\mc{C})$, $\mc{C}(X,Y)$ is an algebra. Moreover, for any $X,Y,Z\in\ob(\mc{C})$, there are morphisms $\bdt^Z_{X,Y}$, $\vps_X$ and $S_{X,Y}$. They  satisfy certain axioms. If $\mc{C}(X,Y)\neq 0$ for any objects $X,Y$, then $\mc{C}$ is called connected. The details about cogroupoids can be found in Subsection 1.2. Let $H$ be a Hopf algebra and let $A$ be a left $H$-Galois object. Then there exists a connected cogroupoid $\mc{C}$ with two objects $X,Y$ such that $H=\mc{C}(X,X)$ and $A=\mc{C}(X,Y)$. Therefore, the theory of Hopf-Galois objects is
actually equivalent to the theory of connected cogroupoids. The following is the main result of this paper.

\begin{thm}\label{intro}
Let $\mc{C}$ be a connected cogroupoid satisfying that $S_{X,Y}$ is bijective for any $X,Y\in \ob(\mc{C})$. Let $X$ be an object in $\mc{C}$ such that $\mc{C}(X,X)$ is a twisted CY algebra of dimension $d$ with homological integral $\int^l_{\mc{C}(X,X)}=\kk_\xi$,  where $\xi:\mc{C}(X,X)\ra \kk$ is an algebra homomorphism. Then for any $Y\in \ob(\mc{C})$, $\mc{C}(X,Y)$ is twisted CY of dimension $d$ with Nakayama automorphism $\mu$ defined as
$$\mu(x)=\xi(x_1^{X,X})S_{Y,X}(S_{X,Y}(x_2^{X,Y})),$$
for any $x\in \mc{C}(X,Y)$.
\end{thm}

The homological integral of a twisted CY algebra will be explained in Subsection 1.4. Sweedler's notation for cogroupoids will be recalled in Subsection 1.2.  The assumption that $S_{X,Y}$ is bijective for any $X,Y\in \ob(\mc{C})$ is to make sure that $S_{Y,X}\circ S_{X,Y}$ is an algebra automorphism. This assumption is not so restrictive, the examples of  cogroupoids mentioned in \cite{bi1} all satisfy this assumption. Actually, if $\mc{C}$ is a connected cogroupoid such that for some object $X$, $\mc{C}(X,X)$ is a Hopf algebra with bijective  antipode, then $S_{X,Y}$ is bijective for any objects $X$, $Y$.  This will be explained in Remark \ref{antipode}. As a direct corollary of the main theorem, we obtain the following description about the CY property of Hopf-Galois objects.

\begin{cor}
Let $H$ be  a twisted CY Hopf algebra of dimension $d$ with bijective antipode. Then any Hopf-Galois object $A$ over $H$ still is a twisted CY algebra of dimension $d$.
\end{cor}

This result generalizes a number of previous results in the literature, such as those of \cite{hvz} for Sridharan algebras, or \cite{yvz}.

An important class of Hopf-Galois objects are cleft Galois objects. Let $H$ be a Hopf algebra.  A left $H$-cleft object is a left $H$-Galois object $A$, such that $A\cong H$ as left $H$-comodules. Left $H$-cleft objects are just the algebras $H_{\sg^{-1}}$, where $\sg$ is a 2-cocycle on $H$ (the definition of the algebras $H_{\sg^{-1}}$ can be found in Subsection 3.1). The algebras $H_{\sg^{-1}}$ are a part of the so-called 2-cocycle cogroupoid. As an application of Theorem \ref{intro}, we obtain that if $H$ is a Hopf algebra with bijective antipode, such that it is twisted CY, then its left cleft Galois objects are all twisted CY. The dual of this  result coincides with Theorem 2.23 in \cite{yvz}. However, the proof in this paper is more direct.

As another application, we obtain that Hopf-Galois objects of the algebras $\mc{B}(E)$, namely the algebras $\mc{B}(E,F)$,  are all twisted CY.
\begin{thm}
Let $E\in \GL_m(\kk)$, $F\in \GL_n(\kk)$ with $m,n\le 2$, and $\tr(E^{-1}E^t)=\tr(F^{-1}F^t)$. The algebra $\mc{B}(E,F)$ is a twisted CY algebra with Nakayama automorphism $\mu$ defined by $$\mu(u)=(E^t)^{-1}Eu(F^t)^{-1}F.$$
\end{thm}

This paper is organized as follows. In Section \ref{s1}, we recall necessary preliminaries and fix some notations. Our main result is proved in  Section \ref{s2}. In the final Section \ref{s3},  we give two applications of the main theorem.

\section{Preliminaries }\label{s1}

Throughout this paper, the unadorned  tensor $\ot$
means $\ot_\kk$ and $\Hom$ means $\Hom_\kk$.

Given an algebra $A$, we write $A^{op}$ for the
opposite algebra of $A$ and $A^e$ for the enveloping algebra $A\ot
A^{op}$. The category of left (resp. right) $A$-modules is denoted by $\Mod$-$A$ (resp. $\Mod$-$A^{op}$). An $A$-bimodule can be identified with an $A^e$-module, that is, an object in $\Mod$-$A^e$

For an $A$-bimodule $M$ and two algebra automorphisms $\mu$ and $\nu$, we let $^\mu M^\nu$ denote the $A$-bimodule such that $^\mu M^\nu\cong M$ as vector spaces, and the bimodule structure is given by
$$a\cdot m \cdot b=\mu(a)m\nu(b),$$
for all $a,b\in A$ and $m\in M$. If one of the automorphisms is the identity, we will omit it. It is well-known that  $A^\mu\cong {}^{\mu^{-1}} A$ as
$A$-bimodules, and that $A^\mu\cong A$ as $A$-bimodules if and only
if $\mu$ is an inner automorphism.


\subsection{Hopf-Galois objects}
We recall the definition of Hopf-Galois objects.
\begin{defn}\label{hg}
Let $H$ be a Hopf algebra. A \textit{left $H$-Galois object} is a left $H$-comodule algebra $A\neq (0)$ such that if $\al:A\ra H\ot A$ denotes the coaction of $H$ on $A$, the linear map
$$\kappa_l: A\ot A\xra{\al\ot 1_A}H\ot A\ot A\xra{1_H\ot m}H\ot A$$
is an isomorphism. A \textit{right $H$-Galois object} is a right $H$-comodule algebra $A\neq (0)$ such that  if $\bt:A\ra A\ot H$ denotes the coaction of $H$ on $A$, the linear map
$$\kappa_r: A\ot A\xra{1_A\ot \bt}A\ot A\ot H\xra{m\ot 1_H}A\ot H$$
is an isomorphism.

If $L$ is another Hopf algebra, an \textit{$H$-$L$-bi-Galois object} is an $H$-$L$-bicomodule algebra which is both a left $H$-Galois object and a right $L$-Galois object.
\end{defn}

\subsection{Cogroupoid}
\begin{defn}A \it{cocategory} $\mc{C}$ consists of:
\begin{itemize}
\item A set of objects $\ob(\mc{C})$.
\item For any $X,Y\in \ob(\mc{C})$, an algebra $\mc{C}(X,Y)$.
\item For any $X,Y,Z\in \ob(\mc{C})$, algebra homomorphisms
$$\bdt^Z_{XY}:\mc{C}(X,Y)\ra \mc{C}(X,Z)\ot \mc{C}(Z,Y) \t{ and } \vps_X:\mc{C}(X,X)\ra \kk$$
such that for any $X,Y,Z,T\in \ob(\mc{C})$, the following diagrams commute:
$$\begin{CD}\mc{C}(X,Y)@>\bdt^Z_{X,Y}>>\mc{C}(X,Z)\ot \mc{C}(Z,Y)\\
@V\bdt^T_{X,Y}VV@V\bdt^T_{X,Z}\ot 1VV\\
\mc{C}(X,T)\ot \mc{C}(T,Y)@>1\ot \bdt^Z_{T,Y}>>\mc{C}(X,T)\ot \mc{C}(T,Z)\ot \mc{C}(Z,Y)\end{CD}$$
\xymatrix
{\mc{C}(X,Y)\ar@{=}[rd]\ar^{\bdt^Y_{X,Y}}[d]&\\
\mc{C}(X,Y)\ot\mc{C}(Y,Y)\ar^-{1\ot \vps_Y}[r]&\mc{C}(X,Y)}\hspace{2mm}
\xymatrix
{\mc{C}(X,Y)\ar@{=}[rd]\ar^{\bdt^X_{X,Y}}[d]&\\
\mc{C}(X,X)\ot\mc{C}(X,Y)\ar^-{\vps_X \ot 1 }[r]&\mc{C}(X,Y).}
\end{itemize}
\end{defn}
Thus a cocategory with one object is just a bialgebra.

A cocategory $\mc{C}$ is said to be \it{connected} if $\mc{C}(X,Y)$ is a non zero algebra for any $X,Y\in \ob(\mc{C})$.

\begin{defn}\label{defn cogoupoid}A \it{cogroupoid} $\mc{C}$ consists of a cocategory $\mc{C}$ together with, for any $X,Y\in \ob(\mc{C})$, linear maps
$$S_{X,Y}:\mc{C}(X,Y)\longrightarrow \mc{C}(Y,X)$$
such that for any $X,Y\in \mc{C}$, the following diagrams commute:
$$\xymatrix{\mc{C}(X,X)\ar[d]_{\bdt_{X,X}^Y}\ar[r]^-{\vps_X}&\kk\ar[r]^-u&\mc{C}(X,Y)\\
\mc{C}(X,Y)\ot\mc{C}(Y,X)\ar[rr]^{1\ot S_{Y,X}}&&\mc{C}(X,Y)\ot\mc{C}(X,Y)\ar[u]^m}$$
$$\xymatrix{\mc{C}(X,X)\ar[d]_{\bdt_{X,X}^Y}\ar[r]^-{\vps_X}&\kk\ar[r]^-u&\mc{C}(Y,X)\\
\mc{C}(X,Y)\ot\mc{C}(Y,X)\ar[rr]^{S_{X,Y}\ot 1}&&\mc{C}(Y,X)\ot\mc{C}(Y,X)\ar[u]^m}$$
\end{defn}
From the definition, we can see that $\mc{C}(X,X)$ is a Hopf algebra for each object $X\in \mc{C}$.

We use Sweedler's notation for cogroupoids. Let $\mc{C}$ be a cogroupoid. For $a^{X,Y}\in \mc{C}(X,Y)$, we write
$$\bdt^Z_{X,Y}(a^{X,Y})=a^{X,Z}_1\ot a^{Z,Y}_2.$$
Now the cogroupoid axioms are
$$(\bdt_{X,Z}^T\ot 1)\circ \bdt^Z_{X,Y}(a^{X,Y})=a^{X,T}_1\ot a^{T,Z}_2\ot a^{Z,Y}_3=(1\ot \bdt_{T,Y}^Z)\circ \bdt^Z_{X,Y}(a^{X,Y});$$
$$\vps_X(a_1^{X,X})a_2^{X,Y}=a^{X,Y}=a_1^{X,Y}\vps_Y(a_2^{Y,Y});$$
$$S_{X,Y}(a_1^{X,Y})a_2^{Y,X}=\vps_X(a_1^{X,X})1=a_1^{X,Y}S_{Y,X}(a_2^{Y,X}).$$

The following is Proposition 2.13 in \cite{bi1}. It describes properties of the ``antipodes''.
\begin{lem}
Let $\mc{C}$ be a cogroupoid and let $X,Y\in\ob(\mc{C})$.
\begin{enumerate}
\item $S_{Y,X}:\mc{C}(Y,X)\ra \mc{C}(X,Y)^{op}$ is an algebra homomorphism.
\item For any $Z\in \ob(\mc{C})$ and $a^{Y,X}\in\mc{C}(Y,X)$,
$$\bdt_{X,Y}^Z(S_{Y,X}(a^{Y,X}))=S_{Z,X}(a_2^{Z,X})\ot S_{Y,Z}(a_1^{Y,X}).$$

\end{enumerate}
\end{lem}

For other basic properties of cogroupoids, we refer to \cite{bi1}.

The following theorem shows the relation between  Hopf-Galois objects and cogroupoids.
\begin{thm}\label{ga}
Let $H$ be a Hopf algebra and let $A$ be a left $H$-Galois object. Then there exists a connected cogroupoid $\mc{C}$ with two objects $X,Y$ such that $H=\mc{C}(X,X)$ and $A=\mc{C}(X,Y)$.
\end{thm}
\proof This is Theorem 2.11 in \cite{bi1} (see also \cite{gr}).

Thus the theory of Hopf-Galois objects is actually equivalent to the theory of connected cogroupoids. In what follows, without otherwise stated, we assume that the cogroupoids mentioned are \textit{connected}.

\subsection{Calabi-Yau algebras}
\begin{defn}\label{tcy}An algebra $A$ is called  a \it{twisted Calabi-Yau algebra  of dimension
$d$} if
\begin{enumerate}\item[(i)] $A$ is \it{homologically smooth}, that is, $A$ has
a bounded resolution by finitely generated projective
$A^e$-modules; \item[(ii)] There is an automorphism $\mu$ of $A$ such that
\begin{equation}\label{cy1}\Ext_{A^e}^i(A,A^e)\cong\begin{cases}0,& i\neq d
\\A^\mu,&i=d\end{cases}\end{equation}
as $A^e$-modules.
\end{enumerate}
If such an automorphism $\mu$ exists, it is unique up to an inner automorphism and is called the \it{Nakayama automorphism} of $A$. A \it{Calabi-Yau algebra} is a twisted Calabi-Yau algebra whose Nakayama automorphism is an inner automorphism. In what follows, Calabi-Yau is abbreviated to CY for short.
\end{defn}

CY algebras were first introduced by Ginzburg in \cite{g2}. Twisted CY algebras include CY algebras as a subclass.  However, the twisted CY property of non-commutative algebras has been studied under other names for many years, even before the definitions of CY algebras. Rigid dualizing complexes of non-commutative algebras were studied in \cite{vdb}. The twisted CY property
was called ``rigid Gorenstein'' in \cite{bz} and was called ``skew Calabi-Yau'' in a recent paper \cite{rrz}.

\subsection{Homological integral}
A standard tool in the theory of the finite dimensional Hopf algebras is the notion
of integrals. Homological integrals for an AS-Gorenstein Hopf algebra introduced in \cite{lwz} is a generalization of this notion. The concept was further extended to a general AS-Gorenstein algebra in \cite{bz}.

The following lemma is probably well-known. Its proof can be found, for example, in \cite[Lemma 2.15]{yvz}.
\begin{lem}\label{kk}
Let $A$ be an augmented algebra such that $A$ is a twisted CY algebra of
dimension $d$ with Nakayama automorphism $\mu$. Then $A$ is of global dimension $d$. Moreover, there is an isomorphism of right $A$-modules
$$\Ext^i_A({}_A\kk,{}_A A)\cong \begin{cases}0,&i\neq d;\\
\kk_\xi,&i=d,\end{cases}$$
where $\xi:A\ra \kk$ is the homomorphism defined by $\xi(h)=\vps(\mu(h))$ for any $h\in A$.
\end{lem}

In a similar way, the following isomorphisms of left $A$-modules holds:
$$\Ext^i_A(\kk_A, A_A)\cong \begin{cases}0,&i\neq d;\\
{}_\eta \kk,&i=d,\end{cases}$$
where $\eta:A\ra \kk$ is the homomorphism defined by $\eta=\vps\circ\mu^{-1}$.

Therefore, if $A$ is a twisted CY augmented algebra, then $A$ has finite global dimension and satisfy the AS-Gorenstein condition. However, $A$ is not necessarily Noetherian.  It is not AS-regular in the sense of \cite[Definition. 1.2]{bz}. We still call the 1-dimensional right $A$-module $\Ext^i_A({}_A\kk,{}_A A)$ right homological integral of $A$ and denoted it by $\int^l_A$.  Similarly, we call the 1-dimensional left $A$-module  $\Ext^d_A(\kk_A, A_A)$  right homological integral of $A$ and denoted it by $\int^r_A$.

\section{The CY property of Hopf-Galois objects}\label{s2}

The aim of this paper is to discuss the CY property of Hopf-Galois objects. Because of the relation between Hopf-Galois objects and connected cogroupoids, this is equivalent to discuss the CY property of connected cogroupoids.  We luckily find that the techniques in \cite{bz} can be used.

Let $\mc{C}$ be a cogroupoid and $X,Y\in \ob(\mc{C})$. Both the morphisms $\bdt_{X,X}^Y:\mc{C}(X,X)\ra \mc{C}(X,Y)\ot \mc{C}(Y,X)$   and $S_{Y,X}:\mc{C}(Y,X)\ra \mc{C}(X,Y)^{op}$ are algebra homomorphisms, so
\begin{equation}\label{mod}(\id\ot S_{Y,X})\circ(\bdt^Y_{X,X}):\mc{C}(X,X)\ra \mc{C}(X,Y)^e(=\mc{C}(X,Y) \ot \mc{C}(X,Y)^{op})\end{equation}
is an algebra homomorphism. This induces a functor $\mc{L}:\Mod$-$\mc{C}(X,Y)^e\ra \Mod$-$\mc{C}(X,X)$.  Let $M$ be a $\mc{C}(X,Y)$-bimodule, that is, an object in $\Mod$-$\mc{C}(X,Y)^e$. The $\mc{C}(X,X)$-module structure of $\mc{L}(M)$ is given by  $$x\ra m=x^{X,Y}_1mS_{Y,X}(x^{Y,X}_2),$$
for any $m\in M$ and $x\in \mc{C}(X,X)$.

In this paper, we view $\mc{C}(X,Y)^e$ as a left and right $\mc{C}(X,Y)^e$-module respectively in the following ways:
\begin{equation}(a\ot b)\ra (x\ot y)=a\ced x\ot y\ced b,\end{equation}
and \begin{equation}\label{lm}(x\ot y)\la(a\ot b) =x\ced a\ot b\ced y.\end{equation}
for any $x\ot y$, $a\ot b \in \mc{C}(X,Y)^e$. Then  $\mc{L}(\mc{C}(X,Y)^e)$ is a left $\mc{C}(X,X)$-module with the module structure given by
$$a\ra (x\ot y)=a^{X,Y}_1\ced x\ot y\ced S_{Y,X}(a^{Y,X}_2)$$
for any $a\in \mc{C}(X,X)$ and $x\ot y\in \mc{C}(X,Y)$. $\mc{L}(\mc{C}(X,Y)^e)$ is actually a $\mc{C}(X,X)$-$\mc{C}(X,Y)^e$-bimodule.

Moreover, the right $\mc{C}(X,Y)^e$-module structure defined in (\ref{lm}) and the algebra homomorphism (\ref{mod}) make $\mc{C}(X,Y)^e$ a right $\mc{C}(X,X)$-module. That is,
\begin{equation}\label{left}(x\ot y)\la a=x\ced a^{X,Y}_1\ot S_{Y,X}(a^{Y,X}_2)\ced y\end{equation}
for any $x\in \mc{C}(X,X)$ and $a\ot b\in \mc{C}(X,Y)^e$.

\begin{lem}\label{R}
\begin{enumerate}
\item[(i)] $\mc{L}$ is an exact functor.
\item[(ii)] $\mc{L}(\mc{C}(X,Y)^e)$ is a free left (right) $\mc{C}(X,X)$-module.
\item[(iii)] $\mc{L}$ preserves injectivity.
\end{enumerate}
\end{lem}
\proof (i) This is clear.

(ii) We define the following morphism:
$$\begin{array}{rcl}
\phi: \mc{C}(X,X)\ot \mc{C}(X,Y) &\lra &\mc{L}(\mc{C}(X,Y)^e)\\
x\ot y &\longmapsto  & x^{X,Y}_1 \ot y \ced S_{Y,X}(x^{Y,X}_2).
\end{array}$$
The morphism $\phi$  is a left $\mc{C}(X,X)$-module homomorphism, where the left  $\mc{C}(X,X)$-module structure of $\mc{C}(X,X)\ot \mc{C}(X,Y)$ is induced from the left module structure of the factor $\mc{C}(X,X)$. We check it in detail as follows:
$$\begin{array}{ccl}
\phi(zx\ot y)&=&z^{X,Y}_1\ced x^{X,Y}_1 \ot y \ced S_{Y,X}(z^{Y,X}_2\ced x^{Y,X}_2 )\\
&=&z^{X,Y}_1\ced x^{X,Y}_1 \ot y \ced S_{Y,X}(x^{Y,X}_2)\ced S_{Y,X}(z^{Y,X}_2)\\
&=&z\ra (x^{X,Y}_1 \ot y \ced S_{Y,X}(x^{Y,X}_2))\\
&=& z\ra \phi (x\ot y).
\end{array}$$
The first and the second equation hold since $\bdt_{X,X}^Y:\mc{C}(X,X)\ra \mc{C}(X,Y)\ot \mc{C}(Y,X)$ and $S_{Y,X}:\mc{C}(Y,X)\ra \mc{C}(X,Y)^{op}$  are algebra homomorphisms.
Moreover, $\phi$ is an isomorphism with inverse
$$\begin{array}{ccl}\mc{L}(\mc{C}(X,Y)^e) &\lra &\mc{C}(X,X)\ot \mc{C}(X,Y)\\x\ot y &\longmapsto& x_1^{X,X} \ot y \ced S_{Y,X}(S_{X,Y}(x_2^{X,Y})).\end{array}$$
Therefore, $\mc{L}(\mc{C}(X,Y)^e)$ is isomorphic to $\mc{C}(X,X)\ot \mc{C}(X,Y)$, it is a free left $\mc{C}(X,X)$-module.

Similarly, we have the following isomorphism of right $\mc{C}(X,X)$-modules:
$$\begin{array}{rcl}
\varphi: \mc{C}(X,X)\ot \mc{C}(X,Y) &\lra &\mc{L}(\mc{C}(X,Y)^e)\\
x\ot y &\longmapsto  & x^{X,Y}_1 \ot S_{Y,X}(x^{Y,X}_2) \ced y.
\end{array}$$
The right  $\mc{C}(X,X)$-module structure of $\mc{C}(X,X)\ot \mc{C}(X,Y)$ is induced from the one of the factor $\mc{C}(X,X)$. The inverse of $\varphi$ is given by
$$\begin{array}{ccl}\mc{L}(\mc{C}(X,Y)^e) &\lra &\mc{C}(X,X)\ot \mc{C}(X,Y)\\x\ot y &\longmapsto& x_1^{X,X} \ot x_2^{X,Y}\ced y.\end{array}$$
Therefore, the module $\mc{L}(\mc{C}(X,Y)^e)$ is also a free right $\mc{C}(X,X)$-module.

(iii) We first show that $\Hom_{\mc{C}(X,Y)^e}({}_{\mc{C}(X,Y)^e}\mc{C}(X,Y)^e_{\mc{C}(X,X)},M)\cong \mc{L}(M)$ as left $\mc{C}(X,X)$-modules for any $\mc{C}(X,Y)$-bimodule $M$. Let $\psi$ be a morphism defined as
$$\begin{array}{rcl}
\psi: \Hom_{\mc{C}(X,Y)^e}({}_{\mc{C}(X,Y)^e}(\mc{C}(X,Y)^e)_{\mc{C}(X,X)},M) &\lra & \mc{L}(M)\\
f &\longmapsto  & f(1\ot 1).
\end{array}$$
We only need to show that $\psi$ is a left $\mc{C}(X,X)$-module homomorphism. For any $x\in \mc{C}(X,X)$, the following equations hold:
$$\begin{array}{ccl}\psi(x\cdot f)&=&(x\cdot f)(1\ot 1)=f((1\ot 1)\la x)\\
&=&f(x^{X,Y}_1\ot S_{Y,X}(x^{Y,X}_2))=x^{X,Y}_1\ced f(1\ot 1)\ced S_{Y,X}(x^{Y,X}_2)\\
&=&x\ra f(1\ot 1) =x\ra \psi(f).
\end{array}$$

Now let $M$ be an injective left $\mc{C}(X,Y)^e$-module. We have $$\begin{array}{cl}
&\Hom_{\mc{C}(X,X)}(-,\mc{L}(M))\\
\cong& \Hom_{\mc{C}(X,X)}(-,\Hom_{\mc{C}(X,Y)^e}({}_{\mc{C}(X,Y)^e}(\mc{C}(X,Y)^e)_{\mc{C}(X,X)},M))\\
\cong &\Hom_{\mc{C}(X,Y)^e}({}_{\mc{C}(X,Y)^e}(\mc{C}(X,Y)^e)_{\mc{C}(X,X)}\ot_{\mc{C}(X,X)}-,M).
\end{array}$$
The last isomorphism follows from the Hom-$\ot$ adjunction. $M$ is an injective $\mc{C}(X,Y)^e$-module and ${}_{\mc{C}(X,Y)^e}(\mc{C}(X,Y)^e)_{\mc{C}(X,X)}$ is a free right $\mc{C}(X,X)$-module by (ii).  So the functor $\Hom_{\mc{C}(X,Y)^e}({}_{\mc{C}(X,Y)^e}(\mc{C}(X,Y)^e)_{\mc{C}(X,X)}\ot_{\mc{C}(X,X)}-,M)$ is exact. Consequently, the functor $\Hom_{\mc{C}(X,X)}(-,\mc{L}(M)$ is exact, that is,  $\mc{L}(M)$ is an injective $\mc{C}(X,X)$-module.
\qed

The following is the key lemma of this paper.

\begin{lem}\label{kk}Let $\mc{C}$ be a cogroupoid, and $X,Y\in \ob(\mc{C})$. The algebra $\mc{C}(X,X)$ is a Hopf algebra, then we have the trivial $\mc{C}(X,X)$-module $\kk$. Let $M$ be a $\mc{C}(X,Y)$-bimodule.
\begin{enumerate}
\item[(i)]$\Hom_{\mc{C}(X,X)}(\kk,\mc{L}(M))= \Hom_{\mc{C}(X,Y)^e}(\mc{C}(X,Y),M).$
\item[(ii)] $\Ext^i_{\mc{C}(X,Y)^e}(\mc{C}(X,Y),M)\cong \Ext^i_{\mc{C}(X,X)}(\kk,\mc{L}(M)),$
for all $i\le 0$.
\end{enumerate}
\end{lem}
\proof This lemma can be obtained by combining \cite[Theorem 7.1]{bi1} and \cite[Lemma 2.4]{bz}(\cite{gk}). Here we give its proof  for the sake of completeness.

(i) By definition,
$$\Hom_{\mc{C}(X,X)}(\kk,\mc{L}(M))=\{m\in M| x\ra m=\vps(x)m, \t{ for any } x\in \mc{C}(X,X)\}$$
and
$$\Hom_{\mc{C}(X,Y)^e}(\mc{C}(X,Y),M)=\{m\in M|ym=my, \t{ for any } y\in \mc{C}(X,Y)\}.$$

If $ym=my$ for all $y\in \mc{C}(X,Y)$, then for any $x\in \mc{C}(X,X)$,
$$\begin{array}{ccl}x\ra m&=&x^{X,Y}_1m S_{Y,X}(x^{Y,X}_2)\\
&=&x^{X,Y}_1(S_{Y,X}(x^{Y,X}_2)m) \\
&=&(x^{X,Y}_1\ced S_{Y,X}(x^{Y,X}_2))m\\
&=&\vps(x) m.\end{array}$$

Conversely, if $x\ra m=\vps(x)m$ for all $x\in \mc{C}(X,X)$, then for any $y\in \mc{C}(X,Y)$,
$$\begin{array}{ccl}my&=&\vps(y^{X,X}_1)my^{X,Y}_2\\
&=&(y^{X,Y}_1mS_{Y,X}(y^{Y,X}_2))y^{X,Y}_3\\
&=&y^{X,Y}_1m(S_{Y,X}(y^{Y,X}_2)\ced y^{X,Y}_3)\\
&=&y^{X,Y}_1m\vps(y_2^{Y,Y})\\
&=&ym.\end{array}$$
Hence, the proof is complete.

(ii)Let $I_\bullet$ be an injective  $\mc{C}(X,Y)^e$-resolution of $M$. $\Ext^i_{\mc{C}(X,Y)^e}(\mc{C}(X,Y),M)$  are the homologies of the complex $\Hom_{\mc{C}(X,Y)^e}(\mc{C}(X,Y),I_\bullet)$. We know from Lemma \ref{R} that $\mc{L}$ is exact and preserve injectivity. Hence $\mc{L}(I_\bullet)$ is an injective resolution of $\mc{L}(M)$.  So  $\Ext^i_{\mc{C}(X,Y)}(\kk,\mc{L}(M))$ are the homologies of the complex $\Hom_{\mc{C}(X,Y)}(\kk,\mc{L}(I_\bullet))$.  Now, as a consequence of (i), we obtain that
$$\Ext^i_{\mc{C}(X,Y)^e}(\mc{C}(X,Y),M)\cong \Ext^i_{\mc{C}(X,X)}(\kk,\mc{L}(M)),$$
for all $i\le 0$.
\qed

Until the end of the section, we assume any cogroupoid $\mc{C}$ mentioned  satisfies that $S_{X,Y}$ is bijective for any $X,Y\in \ob(\mc{C})$. Under this assumption, the morphism $S_{Y,X}\circ S_{X,Y}$ will be an algebra automorphism of $\mc{C}(X,Y)$.


Let $\mc{C}$ be a cogroupoid. For an object $X\in \mc{C}$, $\mc{C}(X,X)$ is a Hopf algebra.  Let $\xi:\mc{C}(X,X)\ra \kk$ be an algebra homomorphism, then we have the right $\mc{C}(X,X)$-module $\kk_\xi$. Let $Y$ be another object in $\mc{C}$, $\mc{L}(\mc{C}(X,Y)^e)$ is a $\mc{C}(X,X)$-$\mc{C}(X,Y)^e$-bimodule, the right $\mc{C}(X,Y)^e$-module structure is given by the equation (\ref{lm}). Then the tensor product $\kk_\xi \ot_{\mc{C}(X,X)} \mc{L}(\mc{C}(X,Y)^e)$ is a right $\mc{C}(X,Y)^e$-module. We have the following lemma.

\begin{lem}\label{iso} Let $\mc{C}$ be a cogroupoid, and $X,Y\in \ob(\mc{C})$. Let $\xi:\mc{C}(X,X)\ra \kk$ be an algebra homomorphism. Then we have
$$ \kk_\xi \ot_{\mc{C}(X,X)} \mc{L}(\mc{C}(X,Y)^e)\cong  (\mc{C}(X,Y))^\mu$$ as right $\mc{C}(X,Y)^e$-modules, where $\mu$ is the algebra automorphism of $\mc{C}(X,Y)$ defined by $$\mu(x)=\xi(x_1^{X,X})S_{Y,X}(S_{X,Y}(x_2^{X,Y}))$$ for any $x\in \mc{C}(X,Y)$.
\end{lem}
\proof Define a morphism $\vph:\kk_\xi \ot_{\mc{C}(X,X)} \mc{L}(\mc{C}(X,Y)^e)\ra (\mc{C}(X,Y))^\mu$ by $\vph(1\ot x\ot y)=\xi(x_1^{X,X}) y\ced S_{Y,X}(S_{X,Y}(x_2^{X,Y}))$
and a morphism $\psi:(\mc{C}(X,Y))^\mu \ra \kk_\xi \ot_{\mc{C}(X,X)} \mc{L}(\mc{C}(X,Y)^e)$ by $\psi(x)=1\ot 1\ot x$.  We first show that $\vph$ is well-defined. We have
$$\begin{array}{cl}
&\vph(1\ot z\ra (x\ot y))\\
=& \vph(1\ot  (z_1^{X,Y}\ced x\ot y\ced S_{Y,X}(z_2^{Y,X})))\\
=& \xi(z_1^{X,X}\ced x_1^{X,X}) y\ced S_{Y,X}(z_3^{Y,X})\ced S_{Y,X}(S_{X,Y}(z_2^{X,Y}\ced x_2^{X,Y})) \\
=&\xi(z_1^{X,X}\ced x_1^{X,X}) y\ced S_{Y,X}(z_3^{Y,X})\ced S_{Y,X}(S_{X,Y}(z_2^{X,Y}))\ced S_{Y,X}(S_{X,Y}(x_2^{X,Y}))\\
=&\xi(z_1^{X,X}\ced x_1^{X,X}) y\ced S_{Y,X}(S_{X,Y}(z_2^{X,Y})\ced z_3^{Y,X})\ced S_{Y,X}(S_{X,Y}(x_2^{X,Y}))\\
=&\xi(z_1^{X,X}\vps(z_2^{X,X})\ced x_1^{X,X}) y\ced S_{Y,X}(S_{X,Y}(x_2^{X,Y}))\\
=&\xi(z)\xi( x_1^{X,X}) y\ced S_{Y,X}(S_{X,Y}(x_2^{X,Y}))\\
=&\vph(\xi(z)\ot x\ot y).
\end{array}$$
Similar calculations show that $\vph$ and $\psi$ are right $\mc{C}(X,Y)^e$-module
homomorphisms and they are inverse to each other.\qed

\begin{lem}\label{hom}
Let $\mc{C}$ be a cogroupoid, and $X,Y\in \ob(\mc{C})$. If $\mc{C}(X,X)$ is homologically smooth, then so is $\mc{C}(X,Y)$.
\end{lem}
\proof The algebra $\mc{C}(X,X)$ is homologically smooth, then $\mc{C}(X,X)$ admits a bounded resolution of finitely generated projective $\mc{C}(X,X)^e$-modules, say
$$P_d\ra P_{d-1}\ra \cdots \ra P_1\ra P_0\ra \mc{C}(X,X)\ra 0.$$
By tensoring the functor $\kk  \ot_{\mc{C}(X,X)}-$, we get a projective right $\mc{C}(X,X)$-module resolution
$$\kk  \ot_{\mc{C}(X,X)} P_d\ra \kk  \ot_{\mc{C}(X,X)} P_{d-1}\ra \cdots  \ra \kk  \ot_{\mc{C}(X,X)} P_1\ra \kk  \ot_{\mc{C}(X,X)} P_0\ra \kk\ra 0$$
of the trivial $\mc{C}(X,X)$-module $\kk$ (cf. \cite[Lemma 2.4]{bt}). The module  $\mc{L}(\mc{C}(X,Y)^e)$ is $\mc{C}(X,X)$-$\mc{C}(X,Y)^e$-bimodule, and is  free as left $\mc{C}(X,X)$-module by Lemma \ref{R} (ii), hence the functor $-\ot_{\mc{C}(X,X)} \mc{L}(\mc{C}(X,Y)^e)$, from $\Mod$-$\mc{C}(X,X)$ to $\Mod$-$(\mc{C}(X,Y)^e)$ is exact and preserves projectivity. Therefore, we obtain a bounded resolution of projective  right $(\mc{C}(X,Y)^e)$-modules:
$$\kk  \ot_{\mc{C}(X,X)} P_d \ot_{\mc{C}(X,X)} \mc{L}(\mc{C}(X,Y)^e)\ra\cdots  \hspace{55mm} $$
$$\hspace{10mm}\ra \kk  \ot_{\mc{C}(X,X)} P_0\ot_{\mc{C}(X,X)} \mc{L}(\mc{C}(X,Y)^e)\ra \kk\ot_{\mc{C}(X,X)} \mc{L}(\mc{C}(X,Y)^e)\ra 0$$
Moreover, each term $\kk  \ot_{\mc{C}(X,X)} P_i \ot_{\mc{C}(X,X)} \mc{L}(\mc{C}(X,Y)^e)$ is a finitely generated  $\mc{C}(X,Y)^e$-module, since $P_i$ is a finitely generated $\mc{C}(X,X)^e$-module. It follows from Lemma \ref{iso} that  $\kk\ot_{\mc{C}(X,X)} \mc{L}(\mc{C}(X,Y)^e)\cong \mc{C}(X,Y)^\mu$ as right $\mc{C}(X,Y)^e$-modules, where $\mu=S_{Y,X}\circ S_{X,Y}$. In conclusion, we obtain a bounded  resolution of $\mc{C}(X,Y)^\mu$ by finitely generated projective right $\mc{C}(X,Y)^e$-modules:
$$\kk  \ot_{\mc{C}(X,X)} P_d \ot_{\mc{C}(X,X)} \mc{L}(\mc{C}(X,Y)^e)\ra\cdots  \hspace{55mm} $$
$$\hspace{10mm}\ra \kk  \ot_{\mc{C}(X,X)} P_0\ot_{\mc{C}(X,X)} \mc{L}(\mc{C}(X,Y)^e)\ra \mc{C}(X,Y)^\mu\ra 0.$$
For each term $\kk  \ot_{\mc{C}(X,X)} P_i \ot_{\mc{C}(X,X)} \mc{L}(\mc{C}(X,Y)^e)$, it is well-known that if we twist its right $\mc{C}(X,Y)^e$-module structure by an automorphism of $\mc{C}(X,Y)^e$, then it is still a projective right $\mc{C}(X,Y)^e$-module. Therefore, in the above resolution,  if we twist the right module structure of each term by the automorphism $\mu^{-1}\ot 1$, then we obtain a desired bounded resolution of $\mc{C}(X,Y)$ by finitely generated projective  $\mc{C}(X,Y)^e$-modules. The algebra  $\mc{C}(X,Y)$ is homologically smooth.

 \qed

Now we are ready to provide our main result.
\begin{thm}\label{main}
Let $\mc{C}$ be a cogroupoid and let $X\in \ob(\mc{C})$ such that $\mc{C}(X,X)$ is a twisted CY algebra with homological integral $\int^l_{\mc{C}(X,X)}=\kk_\xi$,  where $\xi:\mc{C}(X,X)\ra \kk$ is an algebra homomorphism. Then for any $Y\in \ob(\mc{C})$, $\mc{C}(X,Y)$ is twisted CY with Nakayama automorphism $\mu$ defined as
$$\mu(x)=\xi(x_1^{X,X})S_{Y,X}(S_{X,Y}(x_2^{X,Y})),$$
for any $x\in \mc{C}(X,Y)$.
\end{thm}
\proof By Lemma \ref{kk}, for all $i\le 0$,
$$\Ext^i_{\mc{C}(X,Y)^e}(\mc{C}(X,Y),\mc{C}(X,Y)^e)\cong \Ext^i_{\mc{C}(X,X)}(\kk,\mc{L}(\mc{C}(X,Y)^e)).$$

In the proof of Lemma \ref{R} (ii), we obtain a left $\mc{C}(X,X)$-module isomorphism:
$$\begin{array}{rcl}
\phi: \mc{C}(X,X)\ot \mc{C}(X,Y) &\lra &\mc{L}(\mc{C}(X,Y)^e)\\
x\ot y &\longmapsto  & x^{X,Y}_1 \ot y \ced S_{Y,X}(x^{Y,X}_2).
\end{array}$$
This isomorphism is also an isomorphism of right $\mc{C}(X,Y)^e$-modules if we endow a right $\mc{C}(X,Y)^e$-module structure on $\mc{C}(X,X)\ot \mc{C}(X,Y)$ as follows:
$$(a\ot b)\leftarrow (x\ot y)=a x_1^{X,X}\ot y\ced b\ced S^{Y,X}(S^{X,Y}(x_2^{X,Y})),$$
for any $a\ot b\in \mc{C}(X,X)\ot \mc{C}(X,Y)$ and $x\ot y\in \mc{C}(X,Y)^e$.

Now we obtain the following right $\mc{C}(X,Y)^e$-module isomorphisms:
$$\begin{array}{cl}
&\Ext^i_{\mc{C}(X,X)}(\kk,\mc{L}(\mc{C}(X,Y)^e))\\
\cong& \Ext^i_{\mc{C}(X,X)}(\kk,\mc{C}(X,X)\ot \mc{C}(X,Y))\\
\cong& \Ext^i_{\mc{C}(X,X)}(\kk,\mc{C}(X,X))\ot \mc{C}(X,Y)\\
\cong& \Ext^i_{\mc{C}(X,X)}(\kk,\mc{C}(X,X))\ot_{\mc{C}(X,X)}\mc{C}(X,X)\ot \mc{C}(X,Y)\\
\cong& \Ext^i_{\mc{C}(X,X)}(\kk,\mc{C}(X,X))\ot_{\mc{C}(X,X)} \mc{L}(\mc{C}(X,Y)^e).\end{array}$$

The algebra $\mc{C}(X,X)$ is a twisted CY algebra with homological integral $\int^l_{\mc{C}(X,X)}=\kk_\xi$. Therefore, $\Ext^i_{\mc{C}(X,X)}(\kk, \mc{C}(X,X))=0$ for $i\neq d$ and $\Ext^d_{\mc{C}(X,X)}(\kk, \mc{C}(X,X))=\kk_\xi$. We obtain that, $\Ext^i_{\mc{C}(X,Y)}(\mc{C}(X,Y),\mc{C}(X,Y)^e)=0$ for $i\neq d$ and $$\begin{array}{ccl}\Ext^d_{\mc{C}(X,Y)^e}(\mc{C}(X,Y),(\mc{C}(X,Y))^e)&\cong& \kk_\xi\ot_{\mc{C}(X,X)}\mc{L}(\mc{C}(X,Y)^e)
\\&\cong &\mc{C}(X,Y)^\mu,\end{array}$$
where $\mu$ is the algebra automorphism of $\mc{C}(X,Y)$ defined by $$\mu(x)=\xi(x_1^{X,X})S_{Y,X}(S_{X,Y}(x_2^{X,Y})).$$
The second isomorphism follows from Lemma \ref{iso}.

In conclusion, we have that
$$\Ext^i_{\mc{C}(X,Y)^e}(\mc{C}(X,Y),\mc{C}(X,Y)^e)\cong \begin{cases}0&i\neq d;\\
\mc{C}(X,Y)^\mu&i=d.\end{cases}$$

Moreover, by Lemma \ref{hom}, the algebra $\mc{C}(X,Y)$ is homologically smooth. Now the proof is complete.
\qed

\begin{rem}\label{antipode}
If $\mc{C}$ is a cogroupoid such that $\mc{C}(X,X)$ has bijective antipode for some object $X$, then $S_{X,Y}$ is bijective for any objects $X,Y$ in $\mc{C}$. Indeed, let $\mc{C}$ be a cogroupoid such that $\mc{C}(X,X)$ has bijective antipode for some object $X$, then the tensor category of finite dimensional right comodules over $\mc{C}(X,X)$ has right duals, and so has the category of finite dimensional comodules over $\mc{C}(Y,Y)$ for any object $Y$. Then one can use the right duality to construct inverses to all the antipodes $S_{X,Y}$, similarly to the proof of Proposition 2.23 in \cite{bi1}. Combining Theorem \ref{ga}, we obtain the following direct corollary.
\end{rem}

\begin{cor}
Let $H$ be  a twisted CY Hopf algebra of dimension $d$ with bijective antipode. Then any Hopf-Galois object $A$ over $H$ still is a twisted CY algebra of dimension $d$.
\end{cor}

As a consequence of Theorem \ref{main}, we obtain the following duality between Hochschild homology and cohomology.
\begin{cor}
Let $\mc{C}$ be a cogroupoid and let $X\in \ob(\mc{C})$ such that $\mc{C}(X,X)$ is a twisted CY algebra of dimension $d$ with homological integral $\int^l_H=\kk_\xi$,  where $\xi:\mc{C}(X,X)\ra \kk$ is an algebra homomorphism. Let $Y$ be any object in $\mc{C}$. For any $\mc{C}(X,Y)$-bimodule $M$ and for all $0\se i\se d$,
$$\H^i(\mc{C}(X,Y),M)\cong\H_{d-i}(\mc{C}(X,Y), M^\mu),$$
where
$\mu$ is the algebra automorphism of $\mc{C}(X,Y)$ defined as
$$\mu(x)=\xi(x_1^{X,X})S_{Y,X}(S_{X,Y}(x_2^{X,Y}))$$
for any $x\in \mc{C}(X,Y)$.
\end{cor}
\proof By Theorem \ref{main}, the algebra $\mc{C}(X,Y)$ is twisted CY of dimension $d$. Now this corollary is a  consequence of \cite[Theorem 1]{v1}.

\section{Applications}\label{s3}
\subsection{Cleft Hopf-Galois objects}
An important class of Hopf-Galois objects are cleft Galois objects. In this subsection, we discuss their CY property.  To define cleft Galois objects, we first recall the concept of a 2-cocycle.

Let $H$ be a Hopf algebra. A  \it{(right) 2-cocycle} on $H$ is a convolution invertible linear map $\sg:H\ot H\ra \kk$ satisfying
$$\sg(h_1,k_1)\sg(h_2k_2,l)=\sg(k_1,l_1)\sg(h,k_2l_2);$$
$$\sg(h,1)=\sg(1,h)=\vps(h)$$for all $h,k,l\in H$. The set of 2-cocycles on $H$ is denoted $Z^2(H)$.

 The convolution inverse of $\sg$, denote $\sg^{-1}$, satisfies
$$\sg^{-1}(h_1k_1,l)\sg^{-1}(h_2,k_2)=\sg^{-1}(h,k_1l_1)\sg^{-1}(k_2,l_2)$$
$$\sg^{-1}(h,1)=\sg^{-1}(1,h)=\vps(h)$$ for all $h,k,l\in H$.

Let $H$ be a Hopf algebra and let $\sg$ be a 2-cocycle on $H$. The algebra ${}_\sg H$ is defined to be the vector space $H$ together with the multiplication given by
$$h\centerdot k=\sg(h_1,k_1)h_2k_2,$$
for any $h,k\in H$. ${}_\sg H$  is a right $H$-comodule algebra with $\bdt:{}_\sg H\ra {}_\sg H\ot H$ as a coaction, and is a right $H$-Galois object.

Similarly we have the algebra $H_{\sg^{-1}}$. As a vector space $H_{\sg^{-1}}=H$, its multiplication is given by
$$h\centerdot k=h_1k_1\sg^{-1}(h_2,k_2),$$
for any $h,k\in H$. $ H_{\sg^{-1}}$  is a left $H$-comodule algebra with coaction $\bdt:H_{\sg^{-1}}\ra H\ot H_{\sg^{-1}} $, and $H_{\sg^{-1}}$ is a left $H$-Galois object.

The following theorem can be found in \cite{mon}.

\begin{thm}
Let $H$ be a Hopf algebra and let $A$ be a left $H$-Galois object. The following assertions are equivalent:
\begin{enumerate}
\item[(i)] There exists a 2-cocycle $\sg$ on $H$ such that $A\cong H_{\sg^{-1}}$ as left comodule algebras.
\item[(ii)] $A\cong H$ as left $H$-comodules.
\item[(iii)] There exists a convolution invertible $H$-colinear morphism $\phi: H\ra A$.
\end{enumerate}
\end{thm}
A left $H$-Galois object is said to be \textit{cleft} if it satisfies the above equivalent conditions. There is a similar result for right cleft $H$-Galois objects. Right cleft $H$-Galois objects are just the algebras ${}_\sg H$, where $\sg$ is a 2-cocycle.

The algebras $H_{\sg^{-1}}$ and ${}_\sg H$ are part of the 2-cocycle cogroupoid.

Let $H$ be a Hopf algebra and let $\sg,\tau\in Z^2(H)$. The algebra $H(\sg,\tau)$ is defined to be the vector space $H$ together with the multiplication  given by
\begin{equation}\label{mulco}h\centerdot k=\sg(h_1,k_1)h_2k_2\tau^{-1}(h_3,k_3),\end{equation}
for any $h,k\in H$.

It is easy to see that $H(1,1)$ (Here $1$ stands for $\vps\ot\vps$) is just the algebra $H$. For a 2-cocycle $\sg$, $H(1,\sg)$ is the algebra $H_{\sg^{-1}}$ and $H(\sg,1)$ is the algebra ${}_\sg H$.

Now we recall the necessary structural maps for the 2-cocycle cogroupoid of $H$. For any $\sg,\tau, \om\in Z^2(H)$, define the following maps:
\begin{equation}\label{comulco}\begin{array}{rcl}\bdt_{\sg,\tau}^\om=\bdt:H(\sg,\tau)&\longrightarrow &H(\sg,\om)\ot H(\om,\tau)\\
h&\longmapsto &h_1\ot h_2.
 \end{array}\end{equation}
\begin{equation}\label{couco}\vps_\sg=\vps:H(\sg,\sg)\longrightarrow \kk.\end{equation}
\begin{equation}\label{antico}\begin{array}{rcl}S_{\sg,\tau}:H(\sg,\tau)&\longrightarrow & H(\tau,\sg)\\
h&\longmapsto &\sg(h_1,S(h_2))S(h_3)\tau^{-1}(S(h_4),h_5).\end{array}\end{equation}
If $H$ has bijective antipode, then $S_{\sg,\tau}$ is bijective for any $\sg, \tau\in Z^2(H)$.  The inverse of $S_{\sg,\tau}$ is given as:
\begin{equation}\label{anticoinverse}\begin{array}{rcl}S^{-1}_{\sg,\tau}:H(\tau,\sg)&\longrightarrow & H(\sg,\tau)\\
x&\longmapsto &\sg^{-1}(x_5,S^{-1}(x_4))S^{-1}(x_3)\tau(S^{-1}(x_2),x_1).\end{array}\end{equation}

The \it{$2$-cocycle cogroupoid} of $H$, denoted by $\underline{H}$, is the cogroupoid defined as follows:
\begin{enumerate}
\item[(i)] $\t{ob}(\underline{H})=Z^2(H)$.
\item[(ii)] For $\sg,\tau\in Z^2(H)$, the algebra $\underline{H}(\sg,\tau)$ is the algebra $H(\sg,\tau)$ defined in (\ref{mulco}).
\item[(iii)] The structural maps $\bdt^\bullet_{\bullet,\bullet}$, $\vps_\bullet$ and $S_{\bullet,\bullet}$ are defined in (\ref{comulco}), (\ref{couco}) and (\ref{antico}) respectively.
\end{enumerate}
\cite[Lemma 3.13]{bi1} shows that the maps $\bdt^\bullet_{\bullet,\bullet}$, $\vps_\bullet$ and $S_{\bullet,\bullet}$ indeed satisfy the conditions required for a cogroupoid. It is clear that the 2-cocycle cogroupoid is connected.

As an application of Theorem \ref{main}, we obtain the following theorem.
\begin{thm}\label{app1}
Let $H$ be a Hopf algebra with bijective antipode, such that it is twisted CY  with homological integral $\int^l_H=\kk_\xi$, where $\xi:H\ra \kk$ is an algebra homomorphism. Let $\sg$ be a 2-cocycle on $H$,  then the left cleft Galois object $H_{\sg^{-1}}$ is twisted CY  with Nakayama automorphism $\mu$ defined by
$$\mu(x)=\xi(x^{1,1}_1)S_{\sg,1}(S_{1,\sg}(x_2^{1,\sg})),$$
for any $x\in H_{\sg^{-1}}$.
\end{thm}
\proof The 2-cocycle cogroupoid $\underline{H}$ is connected. Moreover, $S_{\sg,\tau}$ is bijective for any $\sg,\tau\in Z^2(H)$, since $H$ has bijective antipode. The algebra $H$ is just the algebra $H(1,1)$  and the algebra $H_{\sg^{-1}}$ is just the algebra $H(1,\sg)$.  Now this theorem is a direct consequence of Theorem \ref{main}.\qed

\begin{rem}Under the assumption of Theorem \ref{app1}, dually, we obtain that a right cleft Galois object ${}_\sg H$ is twisted CY with Nakayama automorphism $\nu$ defined by $\nu(x)=S^{-1}_{\sg,1}(S^{-1}_{1,\sg}(x^{\sg,1}_1))\xi S(x^{1,1}_2)$. This coincides with Theorem 2.23 in \cite{yvz}.
\end{rem}

\subsection{Quantum groups associated to the bilinear form}
Let $E\in \t{GL}_m(\kk)$. Recall that the algebra $\mc{B}(E)$ \cite{dvl} is the algebra presented by generators $(u_{ij})_{1\se i,j\se m}$ and the relations
$$E^{-1}u^tEu=I_m=uE^{-1}u^tE,$$
where $u$ is the matrix $(u_{ij})_{1\se i,j\se n}$, $u^t$ is the transpose matrix and $I_m$ is the identity matrix. It is a Hopf algebra with the following structures:
$$\bdt(u_{ij})=\sum_{k=1}^n u_{ik}\ot a_{kj},\vps(u_{ij})=\dt_{ij},S(u)=E^{-1}u^tE.$$

The following lemma describes the CY property of the algebras $\mc{B}(E)$.
\begin{lem}\label{B cy}
Let $E\in \GL_m(\kk)$. The algebra $\mc{B}(E)$ is a twisted CY algebra of dimension 3 with Nakayama automorphism $\mu$ defined  by
$$\mu(u)=(E^t)^{-1}Eu(E^t)^{-1}E.$$ Its left homological integral is given by $\int^l_{\mc{B}(E)}=\kk_\xi$, where $\xi$ is the algebra homomorphism defined by $\xi(u)=(E^t)^{-1}E(E^t)^{-1}E$.
\end{lem}
\proof The fact that $\mc{B}(E)$ is a twisted CY algebra was shown in \cite{ww}. Here we show that it can be obtained directly from Theorem 6.1 and Corollary 6.3 in \cite{bi}. By \cite[Theorem 6.1]{bi}, the algebra $\mc{B}(E)$ is homologically smooth and $\Ext^i_{\mc{B}(E)^e}(\mc{B}(E),\mc{B}(E)^e)=0$ for $i\le 4$.

\cite[Corollary 6.3]{bi} shows  that $\Ext^i_{\mc{B}(E)^e}(\mc{B}(E),\mc{B}(E)^e)=\Tor_{3-i}^{\mc{B}(E)^e}(\mc{B}(E),{}^\sg(\mc{B}(E)^e))$ for $i=0,1,2,3$.  where $\sg$ is the algebra automorphism given by $\sg(u)=E^{-1}E^tuE^{-1}E^t$. It is easy to see that ${}^\sg(\mc{B}(E)^e)\cong \mc{B}(E)^e$ as left $\mc{B}(E)^e$-modules. So we have $\Tor_i^{\mc{B}(E)^e}(\mc{B}(E),{}^\sg(\mc{B}(E)^e))=0$ unless $i=0$. However,
$\Tor_0^{\mc{B}(E)^e}(\mc{B}(E),{}^\sg(\mc{B}(E)^e))=\mc{B}(E)\ot_{\mc{B}(E)^e} {}^\sg(\mc{B}(E)^e)\cong \mc{B}(E)^{\sg^{-1}}$.
The automorphism $\sg^{-1}$ is just the automorphism $\mu$ defined in this lemma. In conclusion,
$$\Ext^i_{\mc{B}(E)^e}(\mc{B}(E),\mc{B}(E)^e)=\begin{cases}0&i\neq 3;\\
\mc{B}(E)^\mu&i=3.\end{cases}$$
 Now we conclude that $\mc{B}(E)$ is a twisted CY algebra of dimension 3.

It follows from Lemma \ref{kk} that $\int^l_{\mc{B}(E)}=\kk_\xi$, where $\xi=\vps\circ \mu$. That is,  $\xi(u)=(E^t)^{-1}E(E^t)^{-1}E$.  \qed

\begin{rem}
(i) In \cite{bi}, the author worked over the field $\mathbb{C}$. However, the results in that paper, except those in Subsection 6.4 do not depend on the base field.

(ii) The homological integral of $\mc{B}(E)$ can also be obtained by \cite[Proposition 6.2]{bi}. From that proposition,  we have $\int^r_{\mc{B}(E)}=\Ext^3_{\mc{B}(E)}(\kk_{\mc{B}(E)},\mc{B}(E))\cong {}_\eta\kk,$ where $\eta$ is the algebra homomorphism defined by $\eta(u)=E^{-1}E^tE^{-1}E^t$. Consequently,  $\int^l_{\mc{B}(E)}=\kk_\xi$, where $\xi=\eta\circ S$. That is, $\xi(u)=(E^t)^{-1}E(E^t)^{-1}E$.

\end{rem}

The aim of this subsection is to discuss the CY property of the Hopf-Galois objects of the algebras $\mc{B}(E)$. To describe the Hopf-Galois objects, we recall a cogroupoid $\mc{B}$. The algebras $\mc{B}(E)$ is part of this cogroupoid. Let $E\in \t{GL}_m(\kk)$ and let  $F\in \t{GL}_n(\kk)$. The algebra $\mc{B}(E,F)$ is defined to be  the algebra with generators $u_{ij}$, $1\se i\se m$, $1\se j\se n$, subject to the relations:
\begin{equation}\label{alg}F^{-1}u^t E u=I_n;\;\;uF^{-1}u^t E=I_m.\end{equation}
It is clear that $\mc{B}(E)=\mc{B}(E,E)$.

Now we recall the structural maps for the cogroupoid $\mc{B}$. For any $E\in \t{GL}_m(\kk)$, $F\in \t{GL}_n(\kk)$ and $G\in \t{GL}_p(\kk)$, define the following maps:
\begin{equation}\label{bdt}\begin{array}{rcl}
\bdt_{E,F}^G:\mc{B}(E,F)&\longrightarrow &\mc{B}(E,G)\ot \mc{B}(G,F)\\
u_{ij}&\longmapsto &\sum_{k=1}^p u_{ik}\ot u_{kj},
\end{array}\end{equation}
\begin{equation}\label{bvps}\begin{array}{rcl}
\vps_E:\mc{B}(E)&\longrightarrow &\kk\\
u_{ij}&\longmapsto &\dt_{ij},
\end{array}\end{equation}
\begin{equation}\label{bss}\begin{array}{rcl}
S_{E,F}:\mc{B}(E,F)&\longrightarrow &\mc{B}(F,E)^{op}\\
u&\longmapsto &E^{-1}u^tF.
\end{array}\end{equation}
It is clear that $S_{E,F}$ is bijective.

Lemma 3.2 in \cite{bi1} ensures that with these morphisms we have a cogroupoid. The cogroupoid $\mc{B}$  is defined as follows:
\begin{enumerate}
\item[(i)] $\ob(\mc{B})=\{E\in \t{GL}_m(\kk),m\le1\}$.
\item[(ii)] For $E,F\in \ob(\mc{B})$, the algebra $\mc{B}(E,F)$ is the algebra defined as in (\ref{alg})
\item[(iii)] The structural maps $\bdt^\bullet_{\bullet,\bullet}$, $\vps_\bullet$ and $S_{\bullet,\bullet}$ are defined in (\ref{bdt}), (\ref{bvps}) and (\ref{bss}) respectively.
\end{enumerate}

The cogroupoid $\mc{B}$ is not necessarily connected. However, we have the following lemma.

\begin{lem}\rm{(\cite{bi2},\cite[Lemma 3.4]{bi1})}
Let $E\in \GL_m(\kk)$, $F\in \GL_n(\kk)$ with $m,n\le2$. Then $\mc{B}(E,F)\neq (0)$ if and only if $\tr(E^{-1}E^t)=\tr(F^{-1}F^t)$.
\end{lem}

This lemma induces the following corollary.
\begin{cor}
Let $\lambda\in \kk$. Consider the full subcogroupoid $\mc{B}^\lmd$ of $\mc{B}$ with objects
$$\ob(\mc{B}^\lmd)=\{E\in \GL_n(\kk),m\le2, \tr(E^{-1}E^t)=\lmd\}.$$
Then $\mc{B}^\lmd$ is a connected cogroupoid.
\end{cor}

The classification of the Hopf-Galois objects of the Hopf algebras $\mc{B}(E)$ is given in \cite{au}.

\begin{lem}\label{B}
Let $E\in \GL_m(\kk)$ with $m\le2$.
\begin{enumerate}
\item[(i)] If $F\in \GL_n(\kk)$ satisfies $n\le2$ and $\tr(E^{-1}E^t)=\tr(F^{-1}F^t)$, then $\mc{B}(E,F)$ is a left $\mc{B}(E)$-Galois object;
\item[(ii)] Let $A$ be a left $\mc{B}(E)$-Galois object. Then there exists $F\in \GL_n(\kk)$ with $n\le2$ and $\tr(E^{-1}E^t)=\tr(F^{-1}F^t)$, such that $A\cong \mc{B}(E,F)$;
\item[(iii)] For $F\in \GL_n(\kk)$ and $G\in \GL_p(\kk)$, the algebras $\mc{B}(E,F)$ and $\mc{B}(E,G)$ are isomorphic as left $\mc{B}(E)$-comodule algebras if and only if $n=p$ and there is $P\in \GL_n(\kk)$ such that $F=PGP^t$.

\end{enumerate}
\end{lem}

Now we give the CY property of the algebras $\mc{B}(E,F)$.
\begin{thm}
Let $E\in \GL_m(\kk)$, $F\in \GL_n(\kk)$ with $m,n\le 2$ and $\tr(E^{-1}E^t)=\tr(F^{-1}F^t)$. The algebra $\mc{B}(E,F)$ is a twisted CY algebra with Nakayama automorphism $\mu$ defined by $\mu(u)=(E^t)^{-1}Eu(F^t)^{-1}F$.
\end{thm}
\proof $E$ and $F$ are objects in the connected subcogroupoid $\mc{B}^\lmd$, where $\lmd=\tr(E^{-1}E^t)=\tr(F^{-1}F^t)$. Moreover, for any $E',F'\in \ob(\mc{B}^\lmd)$, the morphism $S_{E',F'}$ is bijective. The Hopf algebra $\mc{B}(E)$ is twisted CY by Lemma \ref{B cy}. Now this theorem is a  consequence of  Theorem \ref{main}.

\subsection*{Acknowledgement} The author sincerely thanks the referee for his/her valuable comments and suggestions that helped
her to improve the paper quite a lot. This work is supported by grants from NSFC (No. 11301126, No.11571316), ZJNSF (No. LQ12A01028, No. LY16A010003).

\vspace{5mm}

\bibliography{}

\end{document}